\documentclass[12pt, reqno]{amsart}
\usepackage{amsmath, amsthm, amscd, amsfonts, amssymb, graphicx, color}
\usepackage[bookmarksnumbered, colorlinks, plainpages]{hyperref}

\newtheorem{theorem}{Theorem}[section]
\newtheorem{lemma}[theorem]{Lemma}

\theoremstyle{definition}

\theoremstyle{remark}
\newtheorem{remark}[theorem]{Remark}
\numberwithin{equation}{section}
\DeclareMathOperator{\im}{Im}

\DeclareMathOperator{\diag}{diag}

\DeclareMathOperator{\tr}{tr}
\DeclareMathOperator{\m-m}{m-m}
\begin{document}
\setcounter{page}{1}

\noindent  \textcolor[rgb]{0.00,0.00,0.99}{www.emis.de/journals/BJMA/}\\
%\noindent  \textcolor[rgb]{0.00,0.00,0.99}{www.emis.de/journals/AFA/}\\[.5in]

\title[On positive definite preserving linear transformations \ldots]{ On positive definite preserving linear transformations of rank $r$ on real symmetric matrices}

\author[D. T. Hieu, H. D. Tuan]{Doan The Hieu$^1$$^{*}$ and -Huynh Dinh Tuan$^2$}

\address{$^{1}$  Hue Geometry Group, College of Education, Hue University;
 34 Le Loi, Hue, Vietnam.}

\email{\textcolor[rgb]{0.00,0.00,0.84}{dthehieu@yahoo.com}}

\address{$^{2}$ Hue Geometry Group, College of Education, Hue University;
 34 Le Loi, Hue, Vietnam.}
\email{\textcolor[rgb]{0.00,0.00,0.84}{galois31416@gmail.com}}

%\dedicatory{This paper is dedicated to Professor ABCD}

\subjclass[2010]{Primary 15A86; Secondary 15A18, 15A04.}

\keywords{Linear preserver problems, Symmetric matrix, Positive definite.}

\date{Received: xxxxxx; Revised: yyyyyy; Accepted: zzzzzz.
\newline \indent $^{*}$ Corresponding author}

\begin{abstract}
We study on what conditions on $B_k,$ \ a linear transformation of rank $r$
\begin{equation}\label{form} T(A)=\sum_{k=1}^r\tr(AB_k)U_k\end{equation}
where $U_k,\ k=1,2,\ldots, r$ are linear independent and all positive definite; is positive definite preserving. We give some first results for this question.
For the case of rank one  and two, the necessary  and sufficient conditions are given. We also give some sufficient conditions for the case of rank $r.$
\end{abstract} \maketitle

\section{Introduction}
One of active topics in linear algebra is the linear preserver problems (LPPs) involving linear transformations on matrix space that have special properties: leaving some functions, subsets, relations \ldots invariant.
 For more details about LPPs: the history, the results and open problems we refer the reader to \cite{horn}, \cite{chi_tsi}, \cite{Li_Pi}, \cite{Loe}, and references therein.

On real symmetric or complex Hermitian matrices, the LPP of positive definiteness is still open and seems to be complicated.
In \cite {tutrhi}, we  solved this problem on real symmetric matrices with some additional assumptions.

In this paper, we consider this problem on real symmetric matrices based on the rank of linear transformations. It is showed that a linear transformation  $T$ of rank $r$ preserving positive definiteness can be expressed in the following form
$$T(A)=\sum_{k=1}^r\tr(AB_k)U_k$$
where $U_k,\ k=1,2,\ldots, r$ are linear independent and all positive definite and $B_k,\ k=1,2,\ldots, r$ are all symmetric.

Of course, any linear transformation of form (\ref{form}) may be not positive definite preserving in general.
We address the question on what conditions on $B_k,$ \ $T$ of form (\ref{form}) is positive definite preserving and give some first results. For the case of rank one  and two, the necessary  and sufficient conditions are given. We also give some sufficient conditions for the case of rank $r.$

%=============================================
\section{Some basic lemmas}
%==========================================

\begin{lemma}\label{lem1}
For any $A\in S_{n}(\mathbb{R})$ of rank $r,$ there exists linear independent (pairwise orthogonal) vectors ${\bf x}_1, {\bf x}_2, \ldots, {\bf x}_r$ such that $A=\sum_{i=1}^{r}k_i{\bf x}_i{\bf x}_i^{t},\ k_i\in \{-1,1\}.$ Moreover if $A$ is positive semi-definite, then $A=\sum_{i=1}^{r}{\bf x}_i{\bf x}_i^t.$
\end{lemma}
%==============================================

\begin{lemma} \label{lem2}
 Let $B\in S_n(\mathbb{R}).$
\begin{enumerate}
\item If $B$ is non-zero positive semi-definite, then for every positive definite matrix $A\in S_n(\mathbb{R}),$ we have $\tr(AB)>0.$
\item If $B$ is not positive semi-definite, then there exists positive definite matrix $A\in S_n(\mathbb{R}),$ such that $\tr(AB)<0.$
\end{enumerate}
\end{lemma}
\begin{proof}
\begin{enumerate}
\item By Lemma \ref{lem1}, $A=\sum_{i=1}^{n}x_ix_i^t$, $B=\sum_{i=1}^{r}y_iy_i^t,$ where $\{x_1, x_2,\ldots,x_n\}$ and $\{y_1,y_2,\ldots,y_r\}$ are  systems of linear independent vectors. Then, $AB=\sum_{i=1}^{n}x_ix_i^t\sum_{i=1}^{r}y_iy_i^t=\sum_{i,j}x_ix_i^ty_jy_j^t$. It is not hard to check that
$$\tr(AB)=\sum_{i,j}\tr(x_ix_i^ty_jy_j^t)=\sum_{i,j}\langle x_i, y_j\rangle^2\ge 0.$$
Since $\{x_1, x_2,\ldots,x_n\}$ is a basis of $\Bbb R^n,$ we have $\tr(AB)>0.$
\item
Suppose $QBQ^t=D=\diag(\mu_1,\mu_2,\ldots,\mu_n),$ where $Q$ is an orthogonal matrix. Since $B$ is not positive semi-definite, there exists a diagonal matrix $C>0,$ such that $\tr(CD)<0$ and let $A=Q^tCQ.$ Obviously, $A>0$ and $\tr(AB)=\tr(DC)<0$.
\end{enumerate}
\end{proof}

%-=========================================
The following lemma is well-known in the literature of the theory of quadratic forms.
\begin{lemma} \label{lem3}
Let $A\in S_{n}(\mathbb{R})$ be positive definite and $B\in S_n(\mathbb{R})$. Then there exists an invertible matrix $W,$ such that $WAW^t=I_n$ and $WBW^t=\diag(\mu_1, \mu_2,\ldots,\mu_n).$
\end{lemma}

\begin{remark}  Since $WA=(W^t)^{-1},\ (W^t)^{-1}(A^{-1}B)W^t=diag(\mu_1, \mu_2,\ldots,\mu_n).$
Thus, $\mu_1, \mu_2,\ldots,\mu_n$ are eigenvalues of $A^{-1}B.$
\end{remark}

%====================================================
\section{Positive definite preserving linear transformations of rank $r$}

Let $T:S_{n}(\mathbb{R})\longrightarrow S_{n}(\mathbb{R})$ be a linear transformation of rank $r$ and $\{U_1, U_2,\ldots, U_r\}$ is a basis of $\im T.$ Suppose $T(E_{ii})=\sum_{k=1}^rb^k_{ii}U_k$ and $T(E_{ij}+E_{ji})=\sum_{k=1}^r(b^k_{ij}+b^k_{ji})U_k,\ b^k_{ij}=b^k_{ji},\ k=1,2,\ldots r;\ i,j=1,2,\cdots,n.$
 Let $B_k=(b^k_{ij})_{n\times n},$ then for every $A\in S_{n}(\mathbb{R})$ we can verify
$$T(A)=\sum_{k=1}^r\tr(AB_k)U_k.$$
Of course, a transformation of form (\ref{form}) is linear. Moreover, if $T$  is positive definite preserving, then $U_k,\ k=1,2,\ldots, r$ can be chosen to be positive definite.
\subsection{The case of rank one and two}
By virtue of Lemma \ref{lem2}, the case of rank 1 is easy to prove.
\begin{theorem}
A linear transformation $T$ of rank 1 is positive definite preserving if and only if for every $A\in S_n(\Bbb R),$ \ $T$ has the form.
\begin{equation}\label{form1}T(A)=\tr(AB)U,\end{equation}
where $U>0$ and $B$ is a non-zero positive semi-definite matrix.
\end{theorem}

 In the rest of this subsection, we give the necessary and sufficient condition for the case $T$ is of rank 2.

Consider a linear transformation of rank 2 on $S_n(\Bbb R)$
\begin{equation}\label{form2}T(A)=\tr(AB_1)U_1+\tr(AB_2)U_2,\end{equation}
where $U_1, U_2$ are linear independent and positive definite. By virtue of Lemma \ref{lem3}, there exists an invertible matrix $W$ such that $WU_1W^t=I_n$ and $WU_2W^t=\diag(\mu_1, \mu_2,\ldots, \mu_n);\ \ \mu_i>0,\ i=1,2,\ldots, n.$ Let $\mu_{\min}=\min\{\mu_1, \mu_2,\ldots, \mu_n\}$ and $\mu_{\max}=\max\{\mu_1, \mu_2,\ldots, \mu_n\},$ then we have

\begin{theorem}\label{rank2}
The linear transformation (\ref{form2}) is positive definite preserving if and only if
$B_1+\mu_{\min}B_2$ and $B_1+\mu_{\max}B_2$ are non-zero and positive semi-definite.
\end{theorem}

\begin{proof} Consider the linear transformation $T_1(A)=WT(A)W^t.$ It is easy to see that

$$T_1(A)=\tr(AB_1)I_n+\tr(AB_2)\diag(\mu_1, \mu_2,\ldots, \mu_n),$$
and $T$ is positive definite preserving if and only if $T_1$ is. Thus, we need only to prove the theorem for $T_1.$ First we observe that
\begin{align} \nonumber
&{\ \ \ \ \ }\tr(AB_1)I_n+\tr(AB_2)\diag(\mu_1, \mu_2,\ldots, \mu_n)>0\\ \nonumber
&\Leftrightarrow \tr(AB_1)+\tr(AB_2)\mu_i>0,\ \forall i=1,2,\ldots, n\\
\label{tam}&\Leftrightarrow \tr(AB_1)+\tr(AB_2)\mu_{\min}>0,\tr(AB_1)+\tr(AB_2)\mu_{\max}>0.
\end{align}
But, the fact that (\ref{tam}) is equivalent to $B_1+\mu_{\min}B_2>0$ and $B_1+\mu_{\max}B_2>0,$ is followed by Lemma \ref{lem2}.
\end{proof}
%=====================================================
\subsection{The case of rank $r$}
With a proof similar to the above, we have the following
\begin{theorem}
Let $T$ be the linear transformation of form (\ref{form1}), where $U_i,\ i=1,2,\ldots, U_r$ are linear independent and positive definite. Furthermore, suppose that there exists an invertible matrix $W,$ such that $WU_kW^t=\diag(\lambda_{1k},\lambda_{2k}\ldots,\lambda_{nk}), \ i=1,2,\ldots r,$ then $T$ is positive definite preserving if and only if
$\sum_{k=1}^r\lambda_{ik}B_k,i=1,\ldots,n$ are all non-zero and positive semi-definite.
\end{theorem}
%================================================
By virtue of Lemma \ref{lem2}, we have
\begin{theorem}
Consider a linear transformation of form (\ref{form}), where $U_k>0.$ If $B_k\ge 0,\ i=1,2,\ldots, r$ and are all  non-zero, then $T$ is positive definite pereserving.
\end{theorem}

%===================================================

For a matrix $X\in S_n(\Bbb R),$ denote by $\lambda_{\min}(X),\lambda_{\max}(X)$ the smallest and largest eigenvalues of $X,$ respectively  while $\lambda_{\m-m}(X)$ can be  $\lambda_{\min}(X)$ or $\lambda_{\max}(X).$
The following is a sufficient condition for a linear  transformation of form (\ref{form}) to be  positive definite preserving.
\begin{theorem}
Let $T$ be the linear transformation of form (\ref{form}), where $U_i,\ i=1,2,\ldots, U_r$ are linear independent and positive definite. If
$\sum_{k=1}^r\lambda_{\m-m}(U_k)B_k$ are all non-zero and positive semi-definite, then $T$ is  positive definite preserving.
\end{theorem}

\begin{proof} By the assumption and Lemma \ref{lem2}, we have $\sum_{k=1}^r\tr(AB_k)\lambda_{\m-m}(U_k)>0$ for any positive definite matrix  $A.$ But,
\begin{align}\lambda_{\min}\left[\sum_{k=1}^r\tr(AB_k)U_k\right]&\ge \sum_{k=1}^r\lambda_{\min}[\tr(AB_k)U_k]\\
&\ge \min\{\sum_{k=1}^r\tr(AB_k)\lambda_{\m-m}(U_k)\}>0.\end{align}
The theorem is proved.
\end{proof}
%%%%%%%%%%%%%%%%%%%%%%%%%%%
%{\bf Acknkowledgements.}

\end{document}